\documentclass[a4paper,11pt,reqno]{amsart}

\usepackage[a4paper,margin=3cm]{geometry}
\usepackage{microtype} 

\usepackage{mathtools,amssymb,amsmath,amsthm}
\usepackage{xfrac}
\usepackage{stmaryrd,mathrsfs,bm,yfonts}
\usepackage{esint}

\usepackage[dvipsnames]{xcolor}
\usepackage[pagebackref,colorlinks=true,linkcolor=MidnightBlue,citecolor=Black,urlcolor=MidnightBlue]{hyperref}
\usepackage[capitalise,nameinlink]{cleveref}
\usepackage{doi}

\usepackage{tikz,tikz-3dplot}

\newtheorem{theorem}{Theorem}[section]
\newtheorem{lemma}[theorem]{Lemma}
\newtheorem{proposition}[theorem]{Proposition}
\newtheorem{corollary}[theorem]{Corollary}

\theoremstyle{definition}
\newtheorem{definition}[theorem]{Definition}
\newtheorem{remark}[theorem]{Remark}

\numberwithin{equation}{section}

\newcommand{\N}{\mathbb{N}}

\newcommand{\R}{\mathbb{R}}

\newcommand{\A}{\mathcal{A}}

\newcommand{\Dir}{\mathrm{Dir}}

\newcommand{\eps}{\varepsilon}

\newcommand{\dist}{\mathrm{dist}}

\renewcommand{\epsilon}{\varepsilon}

\def\XXint#1#2#3{{\setbox0=\hbox{$#1{#2#3}{\int}$ }
  \vcenter{\hbox{$#2#3$ }}\kern-.6\wd0}}

\DeclareMathOperator{\Div}{div}

\newcommand{\Iq}{{\mathcal{A}}_Q}
\def\a#1{\left\llbracket{#1}\right\rrbracket}

\newcommand\cH{{\mathcal{H}}}

\newcommand{\de}{\partial}

\makeatletter
\renewcommand{\tocsection}[3]{%
  \indentlabel{\@ifnotempty{#2}{\bfseries\ignorespaces\makebox[\@ifempty{#1}{25pt}{75pt}][l]{#1 #2\quad}}}\bfseries#3}
\renewcommand{\tocsubsection}[3]{%
  \indentlabel{\@ifnotempty{#2}{\ignorespaces\makebox[30pt][l]{#1 #2\quad}}}#3}
\newcommand\@dotsep{4}
\def\@tocline#1#2#3#4#5#6#7{\relax
  \ifnum #1>\c@tocdepth \else
    \par \addpenalty\@secpenalty\addvspace{#2}%
    \begingroup \hyphenpenalty\@M
    \@ifempty{#4}{\@tempdima\csname r@tocindent\number#1\endcsname\relax}{\@tempdima#4\relax}%
    \parindent\z@ \leftskip#3\relax \advance\leftskip\@tempdima\relax
    \rightskip\@pnumwidth plus1em \parfillskip-\@pnumwidth
    #5\leavevmode\hskip-\@tempdima{#6}\nobreak
    \leaders\hbox{$\m@th\mkern \@dotsep mu\hbox{.}\mkern \@dotsep mu$}\hfill
    \nobreak
    \hbox to\@pnumwidth{\@tocpagenum{\ifnum#1=1\bfseries\fi#7}}\par
    \nobreak
    \endgroup
  \fi}
\AtBeginDocument{%
\expandafter\renewcommand\csname r@tocindent0\endcsname{0pt}}
\def\l@section{\@tocline{1}{0pt}{1pc}{25pt}{}}
\def\l@subsection{\@tocline{2}{0pt}{2pc}{30pt}{}}
\makeatother

\title[Carleman for $Q$-valued maps]{Carleman estimates for stationary $Q$-valued maps: a variational approach}

\author[A. Halavati]{Aria Halavati}
\address{Courant Institute of Mathematical Sciences,
New York University.
251 Mercer Street,
New York, NY 10012-1185}
\email{aria.halavati@cims.nyu.edu}

\address{Bocconi University, Department of Decision Sciences, Via Guglielmo R\"ontgen 1, 20136 Milano,
Italy}
\email{aria.halavati@unibocconi.it}
\author[L. Spolaor]{Luca Spolaor}
\address{Department of Mathematics, UC San Diego, AP\&M, La Jolla, California, 92093, USA}
\email{lspolaor@ucsd.edu}

\begin{document}

\maketitle

\begin{abstract}
We prove a Carleman-type estimate for Dirichlet-stationary multivalued functions and apply it to give a different proof of the optimal dimension of the singular set of Dir-minimizing multivalued functions, originally due to Almgren and to De Lellis–Spadaro.
\end{abstract}

\tableofcontents

\section{Introduction}

Recall that a $Q$-valued map $f\in W^{1,2}(\Omega,\Iq(\R^m))$, $\Omega\subset \R^n$ open, is \emph{Dir-stationary (with respect to outer and inner variations of the Dirichlet energy)} if it is a critical point of the Dirichlet energy, that is it satisfies an \emph{outer variation formula}
\begin{equation}\label{eq:outer}
    \mathcal O(f, \psi):=\int \sum_i\left[\langle Df_i(x)\,:\, D_x \psi(x,f_i(x))\rangle + \langle Df_i(x)\,:\, D_u \psi(x,f_i(x))\cdot Df_i(x)\rangle\right] \,dx=0\,,
\end{equation}
for every $\psi(x,u) \in C^\infty(\Omega\times \R^m;\R^m)$ with \emph{compact support in $x$} and
\[
|D_u\psi|\leq C<\infty\qquad \text{and}\qquad |\psi|+|D_x\psi|\leq C\, (1+|u|)\,,
\]
and an \emph{inner variation formula}
\begin{equation}\label{eq:inner}
    \mathcal I(f, \phi):=2 \int \sum_{i=1}^Q \langle Df_i \,:\, Df_i\cdot D\phi \rangle - \int |Df|^2 \, \Div \phi = 0\,,\qquad \forall \phi \in C^\infty_c(\Omega, \R^n)\,.
\end{equation}
We will say that a $Q$-valued map $f\in W^{1,2}(\Omega,\Iq(\R^m))$ is \emph{weakly stationary} if $\mathcal O(f, \cdot)=0$, i.e.\ stationary with respect to outer variations only.

For the derivation of these formulas as the Euler–Lagrange equations of the Dirichlet energy and the related Sobolev theory of $Q$-valued maps, see \cite[Sections 2 and 3]{DS}, whose notations we follow.

Multivalued maps that minimize an appropriate Dirichlet energy were introduced by Almgren in \cite{Alm} in his celebrated proof of the optimal bound on the dimension of the singular set of area minimizing currents in high codimension, as the appropriate linearized problem. More recently, De Lellis and Spadaro revisited this theory with metric techniques in \cite{DS,DLS_Currents,DLS_Lp,DLS-center-manifold,DLS_Blowup} (see also \cite{Hi} for minimizers taking values in a smooth compact Riemannian manifold, and \cite{BDPW} for minimizers of the $p$-Dirichlet energy). 

The main difficulty in proving such optimal bound is the presence of branch points: points at which the blow-up is regular but the presence of multiplicity causes singularity to appear. Almgren's innovative insight was to understand how such branch points can be studied with the same techniques that are used in the study of unique continuation properties for elliptic PDEs. However multivalued functions do not satisfy a PDE in the usual sense, and so he had to find a variational approach to unique continuation which gave birth to the so-called frequency function. However, another technique that has been extensively used in the context of standard PDEs to study unique continuation type question is the so called Carleman estimates technique (see for instance \cite{Car39,Aronszajn1957,JK85,KT01,Isa04,Rul15}). In this paper we give a variational proof of such an estimate for $Q$-valued functions and we use it to recover Almgren's optimal bound on their singular set. We point out that Carleman estimates have been used in a similar setting for J-holomorphic maps by Riviere-Tian (see \cite{RT-J-holo}) where they take advantage of the complex structure to turn the problem into a first order elliptic system, our approach is different (partially inspired by the proofs in \cite{H1}).

We remark that, although our methods differ, the information obtained through Carleman estimates is essentially captured by the frequency function approach. The main purpose of this note is to make the community aware of this technique (mainly the two inner and outer variations needed to find the estimates).

\subsection{Main results}
The main result of the paper is the following Carleman-type estimate.

\begin{theorem}[Carleman estimate]\label{thm:Carleman-proposition}
    Let $f\in W^{1,2}(B_1;\A_Q(\R^m))$ be a Dir-stationary $Q$-valued function, $B_1\subset \R^n$, and let $\tau>0$ and $\eta = \frac{2\tau-n+2}{2}$. Then the following estimate holds
    \begin{align}\label{eq:carlo}
        \int_{B_1} \chi \sum_{i=1}^Q \left(\epsilon^2 \frac{|f_i|^2}{|x|^{2\tau + 2 - \epsilon}}+ \frac{|Df_i\cdot x - \eta f_i|^2}{|x|^{2\tau+2}} \right) \leq C\int_{B_1} |D\chi|\sum_{i=1}^Q \left(\frac{|Df_i|^2}{|x|^{2\tau-1}} + \frac{|f|^2}{|x|^{2\tau+1}}\right)\,.
    \end{align}
    for any compactly supported function $\chi \in C^{\infty}_c\left(B_1\setminus \{0\}\right)$. 
\end{theorem}

We remark that in the literature the name Carleman estimate is usually associated to estimates of the form
\begin{equation}\label{eq:carlo1}
\left\||x|^{-\tau}Du\right\|_{L^2(\Omega)} \leq \left\||x|^{-\tau+1}\Delta u\right\|_{L^2(\Omega)}\,,
\end{equation}
(for large $\tau>0$) which are then used to derive expressions of the form \cref{eq:carlo}. However in our case the Laplacian is replaced by inner and outer variations for the Dirichlet energy, so \cref{eq:carlo1} doesn't make sense in our setting and we have to give a variational proof of \cref{eq:carlo}, that is by testing inner and outer variations with proper vector fields.. 

A straightforward consequence of \cref{thm:Carleman-proposition} is the following strong unique continuation result

\begin{theorem}[Strong unique continuation]\label{thm:strong}
    Let $f\in W^{1,2}(B_1, \Iq(\R^m))$, with $B_1 \subset \R^n$, be a Dir-stationary map and suppose that 
    \[
    \lim_{r\to 0} \frac1{r^N} \int_{B_r}|f|^2=0\qquad \forall N\in \N\,.
    \]
    Then $f\equiv Q\a{0}$ in $B_1$.
\end{theorem}

Moreover, we are also able to recover the optimal bound on the singular set of Dir-minimizing multivalued functions. We recall that a point \emph{$x\in \Omega$ is regular} if there exists a neighborhood $B\subset \Omega$ of $x$ and $Q$ analytic functions $f_i \colon B \to \R^n$ such that 
\[
f(y) = \sum_{i=1}^Q \a{f_i(y)}\qquad \text{ for
almost every }y \in B\,,
\]
and either $f_i(x)\neq f_j(x)$ for every $x\in B$ or $f_i \equiv f_j$. The \emph{singular set $\Sigma_f$ of $f$} is the complement in $\Omega$ of the set of regular points.

\begin{theorem}[Dimension of the singular set]\label{thm:dimension} Let $f\in W^{1,2}(\Omega, \Iq(\R^m))$, with $\Omega \subset \R^n$, be a Dir-minimizing map. Then $\dim_{\cH}(\Sigma_f)\leq n-2$. Moreover if $n=2$, then $\Sigma_f$ is locally finite.
\end{theorem}

We remark that \cref{thm:dimension} is already known, see \cite{Alm, DS}. Our main contribution is the use of \cref{eq:carlo} in place of the monotonicity formula for the frequency function to prove them.

Finally, as an instructive remark, we note that by adapting ideas from \cite{HS}, one can reprove the result therein for two-dimensional Dir-stationary multivalued functions, replacing frequency function techniques with our Carleman estimate in combination with Weiss’ energy.

\medskip

\noindent \textbf{Acknowledgments.} The authors would like to thank Camillo De Lellis, Guido De Philippis and Robert V.\ Kohn for their interest in this work and related discussions. L.S. acknowledges the support of the NSF Career Grant DMS 2044954 and A.H. acknowledges the support of NSF grant DMS-2055686 and the Simons Foundation and the European Union: the
European Research Council (ERC), through StG “ANGEVA”, project number: 101076411. Views and opinions
expressed are however those of the authors only and do not necessarily reflect those of the European Union
or the European Research Council. Neither the European Union nor the granting authority can be held
responsible for them.

\textbf{Conflict of interest.} On behalf of all authors, the corresponding author states that there is no conflict of interest.

\textbf{Data availability statement.} This article has no associated data.

\section{Proof of Carleman estimate: \cref{thm:Carleman-proposition}}
 For a compactly supported positive test function $\chi \in C^\infty_c(B_1\setminus\{0\})$ we test the outer variation \cref{eq:outer} with the admissible test function $\displaystyle{\psi(x,u) = \chi\,\frac{u}{|x|^{2\tau}}}$ to obtain
        \begin{align}\label{outer-Q-1}
        \begin{aligned}
            0 = &\int_{B_1}  \sum_{i=1}^Q \left[\chi\frac{|Df_i|^2}{|x|^{2\tau}} + \frac{\langle D \chi\,\colon\,Df_i\rangle \,f_i}{|x|^{2\tau}}-2\,\tau\,\chi\, \frac{\de_r f_i\, f_i}{|x|^{2\tau+1}} \right]\\
            = &\int_{B_1} \sum_{i=1}^Q \left[\chi\frac{|Df_i|^2}{|x|^{2\tau}} + \frac{\langle D \chi\,\colon\,Df_i\rangle \,f_i}{|x|^{2\tau}}+\tau(n-2-2\tau)\chi \frac{|f_i|^2}{|x|^{2\tau+2}}  + \tau \de_r\chi\frac{|f_i|^2}{|x|^{2\tau+1}}\right]\,,
            \end{aligned}
        \end{align}
        where the last equality is obtained by an integration by parts.
        Next we want to test the inner variation \cref{eq:inner} with the admissible vector field $\phi = \chi \,\frac{x}{|x|^{2\tau}}$. We compute
        \begin{gather}
        D \phi=\frac{1}{|x|^{2\tau}}\,D\chi\otimes x+\frac{\chi}{|x|^{2\tau}}\,I-2\tau \, \frac{\chi}{|x|^{2\tau+2}}\,x\otimes x\,,\\
        {\rm div}{(\phi)}=\frac{\de_r\chi}{|x|^{2\tau-1}}+(n-2\tau)\,\frac{\chi}{|x|^{2\tau}}\,.
        \end{gather}
        This implies that
        \begin{align}
        \begin{aligned}\label{inner-Q-1}
            \int_{B_1} \sum_{i=1}^Q \left[\frac{(n-2\tau-2)}{2}\chi\frac{|Df_i|^2}{|x|^{2\tau}} + 2\tau \chi \frac{|\de_r f_i|^2}{|x|^{2\tau}}  + \frac{\partial_r \chi}{2}\frac{|Df_i|^2}{|x|^{2\tau-1}} - \frac{\langle Df_i\,\colon\,D\chi\rangle \, \langle Df_i\,\colon \, x\rangle }{|x|^{2\tau}}\right] = 0\,.
        \end{aligned}
        \end{align}
        Now we name $\eta = \frac{2\tau-n+2}{2}$, multiply \cref{outer-Q-1} by $\eta$ and add it to \cref{inner-Q-1} to see that:
        \begin{align}
        \begin{aligned}\label{pre-carleman}
            \int_{B_1} \chi\sum_{i=1}^{Q}\left[\frac{|\de_r f_i|^2}{|x|^{2\tau}} - \eta^2\frac{|f_i|^2}{|x|^{2\tau+2}}\right] \leq C\left(\frac\eta\tau\right)\int_{B_1} |D\chi|\sum_{i=1}^{Q}\left[\frac{|Df_i|^2}{|x|^{2\tau-1}}+\frac{|f_i|^2}{|x|^{2\tau+1}}\right]\,.
            \end{aligned}
        \end{align}
        This is a first Carleman estimate:  to conclude we need to complete the square on the left hand side. We calculate as follows:
        \begin{align*}
            \int_{B_1} 2\eta^2 \chi \sum_{i=1}^Q \frac{|f_i|^2}{|x|^{2\tau+2}} &= \int_{S^{n-1}}\left(\int_{0}^1 2\eta^2\chi \sum_{i=1}^Q\frac{|f_i|^2}{r^{2\tau-n+3}}\,dr\right)d\theta\notag\\
            &= \int_{S^{n-1}}\left(\int_{0}^1 -\eta\chi \de_r(r^{-2\tau+n+2})\sum_{i=1}^Q|f_i|^2\,dr\right)d\theta\notag\\
            &= \int_{S^{n-1}}\left(\int_{0}^1 2\eta\,\chi\,r^{-2\tau+n-2}\, \sum_{i=1}^Qf_i\de_rf_i\,dr\right)d\theta\notag\\
            &\quad + \int_{S^{n-1}}\left(\int_{0}^1 \eta\,\de_r\chi \,r^{-2\tau+n-2}\sum_{i=1}^Q|f_i|^2\,dr\right)d\theta\,,
        \end{align*}
        that is
        \begin{equation}\label{temp-1}
            \int_{B_1} 2\eta^2 \chi \sum_{i=1}^Q \frac{|f_i|^2}{|x|^{2\tau+2}} -\int_{B_1} 2\eta\,\chi\, \sum_{i=1}^Q\frac{f_i \,Df_i\cdot x}{|x|^{2\tau+2}}=\eta\,\int_{B_1}\,D\chi\cdot x\,\sum_{i=1}^Q\frac{|f_i|^2}{|x|^{2\tau+2}}\,
        \end{equation}
        Then, combining \cref{pre-carleman} and \cref{temp-1}, we see that:
        \begin{align}\label{First-carleman}
            \int_{B_1} \chi\sum_{i=1}^Q\frac{\left|Df_i\cdot x-\eta f_i\right|^2}{|x|^{2\tau+2}} \leq C \int_{B_1} |D\chi| \sum_{i=1}^Q \left(\frac{|Df_i|^2}{|x|^{2\tau-1}} + \frac{|f_i|^2}{|x|^{2\tau+1}}\right)\,.
        \end{align}
        
        To estimate the $L^2$ term on the left hand side of \cref{eq:carlo}, we can see \cref{temp-1} for $\tau-\epsilon$ in place of $\tau$ and the obvious modification for $\eta$:
        \begin{align}\label{epsilon-temp-1}
            \int_{B_1} \chi\sum_{i=1}^Q \left(|\de_r f_i|^2 - (\eta-\epsilon)^2 |f_i|^2\right) \geq -C\int_{B_1} |D\chi| \sum_{i=1}^Q \left(\frac{|Df_i|^2}{|x^{2\tau-1}} + \frac{|f_i|^2}{|x|^{2\tau+1}}\right)\,.
        \end{align}
        Combining this, with \cref{First-carleman} and the following computation, we conclude:
        \begin{align*}
         \int_{B_1} \chi\sum_{i=1}^Q\frac{\left|Df_i.x-\eta f_i\right|^2}{|x|^{2\tau+2-2\epsilon}}
            &\geq \int_{B_1} \chi\sum_{i=1}^Q\frac{|\de_r f_i|^2}{|x|^{2\tau+2-2\epsilon}} - \eta \frac{\de_r |f_i|^2}{|x|^{2\tau+1-2\epsilon}}+ \eta^2\frac{|f_i|^2}{|x|^{2\tau-2\epsilon}}\\
            &\geq \int_{B_1} \chi\sum_{i=1}^Q\frac{|\de_r f_i|^2}{|x|^{2\tau+2-2\epsilon}} + (\eta^2-2\eta(\eta-\epsilon))\frac{|f_i|^2}{|x|^{2\tau-2\epsilon}}\\
            &\stackrel{\cref{epsilon-temp-1}}{\geq} \int_{B_1} \chi\sum_{i=1}^Q\epsilon^2\frac{|f_i|^2}{|x|^{2\tau-2\epsilon}} \\ &\;\;\;\;\;\;\;\;\;-C \int_{B_1} |D\chi| \sum_{i=1}^Q \left(\frac{|Df_i|^2}{|x|^{2\tau-1}} + \frac{|f_i|^2}{|x|^{2\tau+1}}\right)\,.
        \end{align*}
        \qed

\begin{corollary}[Three sphere inequality]\label{three-sphere-multivalued-corollary}
    For any constant $\tau_0>0$, there exists a constant $C(\tau_0)>0$, possibly depending on the dimension, with the following property. For any constant $0< \tau < \tau_0$ (possibly large), any Dir-stationary multifunction $f\in W^{1,2}(B_1;\mathcal{A}_Q)$ and three radii $r_1< r_2 < r_3 < \frac{1-|x|}{2}$ with $\min\left(\frac{r_3}{r_2},\frac{r_2}{r_1}\right)>2$ and $x\in B_1$, the following estimate is true:
    \begin{align*}
        &\left[\frac{1}{1+\log(r_3/r_2)^2} + \frac{1}{1+\log(r_2/r_1)^2}\right]\frac{\|f\|^2_{L^2(B_{2r_2}\setminus B_{r_2}(x))}}{r_2^{2\tau}} \\ &\leq  C(\tau_0)\left( \frac{\|f\|^2_{L^2(B_{2r_1}\setminus B_{r_1}(x))}}{r_1^{2\tau}} + \frac{\|f\|^2_{L^2(B_{2r_3}\setminus B_{r_3}(x))}}{r_3^{2\tau}}\right).
    \end{align*}
    Here $C(\tau_0) \leq e^{K\tau_0}$ for some universal $K>0$.
    \end{corollary}
    
    \begin{proof}
    First testing the outer variation \cref{eq:outer} with $\Psi(x,u) = \phi^2(x)u$ we get the Caccioppoli inequality:
        \begin{align}\label{eq:caccio}
            \int_{B_1} \phi^2|Df|^2 \leq C \int_{B_1}|D\phi|^2|f|^2\,.
        \end{align}
        Without loss of generality we assume that $x$ is the origin. Then we distinguish two cases.
        
        \textit{Case I:} If $\log(r_3/r_2) > \log(r_2/r_1)$ (meaning that $r_2$ is closer to $r_1$ than $r_3$) we rescale such that $r_1$ becomes the unit radius. Then we choose $\chi = 1$ on $B_{1.1 r_3/r_1}\setminus B_{1.9 }$, supported in $B_{1.9r_3\setminus r_1}\setminus B_{1.1}$, and to decay linearly to zero in $B_{1.1}$ and outside $B_{1.1 r_3/r_1}$ so that $|D\chi| \leq C^{-1}$ in $B_{1.9}\setminus B_{1.1 }$ and $|D\chi| \leq C (r_3/r_1)^{-1}$ in $B_{1.9 r_3/r_1}\setminus B_{1.1 r_3/r_1}$. Then we use \cref{thm:Carleman-proposition} and rescale back with $r_1>0$ to see that for any $\epsilon >0$:
        \begin{align*}
            \epsilon^2  \left(\frac{r_2}{r_1}\right)^{2\epsilon}r_2^{-2\tau}\int_{B_{2r_2}\setminus B_{r_2}(x)} |f|^2 \leq C(\tau_0) r_1^{-2\tau}&\int_{B_{1.9 r_1}\setminus B_{1.1 r_1}(x)} |f|^2 + r_1^2|Df|^2 \\ +  C(\tau_0) r_3^{-2\tau}&\int_{B_{1.9 r_3}\setminus B_{1.1 r_3}(x)} |f|^2 + r_3^2|Df|^2\,.
        \end{align*}
        Taking a smooth test function $\phi=1$ on $B_{1.9r_1}\setminus B_{1.1r_1}(x)$ and $\phi = 0$ outside $B_{2r_1}\setminus B_{r_1}(x)$ with $|d\phi|\leq Cr_1^{-1}$ and similarly for $r_3$ and using the Cacciopoli inequality we can bound the gradient terms to get that:
        \begin{align*}
            \epsilon^2  \left(\frac{r_2}{r_1}\right)^{2\epsilon} r_2^{-2\tau}\int_{B_{2r_2}\setminus B_{r_2}(x)} |f|^2 \leq C(\tau_0)\left( r_1^{-2\tau}\int_{B_{2r_1}\setminus B_{r_1}(x)} |f|^2+  r_3^{-2\tau}\int_{B_{2 r_3}\setminus B_{ r_3}(x)} |f|^2\right).
        \end{align*}
        The desired estimate comes from optimizing $\epsilon^2 \left(\frac{r_2}{r_1}\right)^{2\epsilon}$ and plugging in $\epsilon = \frac{1}{\sqrt{\ln(r_2/r_1)^2 + 1}}$.

        \textit{Case II:} If $\log(r_3/r_2) \leq \log(r_2/r_1)$ Then we rescale so that $r_3$ radius becomes unit scale and perform the same analysis. The result follows by adding the two possibilities.
    \end{proof}

\section{The vanishing order and its properties}
This part of our paper is inspired by the work of \cite{KRS-carleman}, and we follow their strategy. We start by defining the vanishing order of a Dir-minimizing multivalued function as follows:

\begin{definition}[Vanishing order]
    Let $f\in W^{1,2}(B_1,\mathcal{A}_Q)$ be a Dir-stationary multi-valued function. Then around any point $x\in B_1$ we define the vanishing degree $\kappa_x$ as follows:
    \begin{align*}
        \kappa_{x,f} = \limsup_{r\to0} \frac{\log\left(\fint_{B_{2r}\setminus B_r(x)} |f|^2 \right)}{2\log(r)}=\lim_{r\to0} \frac{\log\left(\fint_{B_{2r}\setminus B_r(x)} |f|^2 \right)}{2\log(r)}
    \end{align*}
    When it's clear from the context we will drop the subindex $f$.
\end{definition}

\begin{remark}
    We define the vanishing order around all points, however only collapsed points are relevant since non-collapsed points have $0$ as their vanishing order. Also notice that if $f(x)=Q\a{y}$, for some $y\in \R^n$, then a more meaningful definition of vanishing order could be obtained by replacing $|f|^2$ with $|f\ominus y|^2$ in the integral. However, later we will consider function with zero average, and so we will not change the definition here.
    Finally we observe that, if $f\equiv Q \a{0}$ in a neighborhood of a point $x$, then clearly $\kappa_x=\infty$.
\end{remark}

In the next two subsections we will prove that the vanishing order is well defined (i.e., the limit exists) and we will show some of its properties, namely upper semicontinuity and homogeneity of suitable sequences of blow-ups.

\subsection{The vanishing order and strong unique continuation}

In order to prove that the limit above is well defined we will need the following immediate corollary of \cref{thm:Carleman-proposition}.

Using the three sphere inequality first we show that:
\begin{lemma}[Strong unique continuation]\label{lemma:strong-unique-continuation}
    Let $f$ be an average free, Dir-stationary functions on an open and connected domain $\Omega$. If there exists a point $x\in \Omega$ such that $\kappa_x=+\infty$, then $f=Q\a{0}$ in $\Omega$. In Particular \cref{thm:strong} is true.
    \begin{proof}
        Take the set $S = \{x\in\Omega: \kappa_x=\infty\}$ and assume it is non-empty. We aim to show that $S$ is either empty or $S=\Omega$. By contradiction, assume that $\Omega\not=S$. Then for any point $x\in S$ we have:
        \begin{align*}
            \limsup_{r\to0} \frac{\log\left(\fint_{B_{2r}\setminus B_r(x)} |f|^2 \right)}{2\log(r)} = \infty\,.
        \end{align*}
        This means that the $L^2$-norm of $f$ vanishes faster than any polynomial. Hence we can take $r_1 \to 0$ in the Carleman estimate of \cref{thm:Carleman-proposition} to see that for any $\tau>0$ and $r_1 \geq 2r_2$:
        \begin{align*}
            \int_{B_{2r_2}\setminus B_{r_2}(x)} |f|^2 \leq e^{K\tau}(\log(r_2)^2+1) \left(\frac{r_2}{r_1}\right)^{2\tau}\int_{B_{2r_1}\setminus B_{r_1}(x)} |f|^2 \,.
        \end{align*}
        Then if $r_1 \geq e^{2K} r_2$, we can take $\tau\to\infty$ to see that for all $r < \delta_0 \dist(x, \partial \Omega)$ for some fixed $\delta_0>0$:
        \begin{align*}
            \int_{B_{2r}\setminus B_{r}(x)} |f|^2 = 0\,.
        \end{align*}
        Hence we see that $f(x)=Q\llbracket 0\rrbracket$ and $\kappa_x=\infty$ for all $x\in B_{\delta_0\dist(x, \partial \Omega)}(x)$. Hence we have just shown that:
        \begin{align*}
            \text{for all } x\in S \text{ we have } B_{\delta_0\dist(x,\de\Omega)}(x)\subset S\,.
        \end{align*}
        Since $\Omega$ is connected, we can find a continuous path $\gamma:[0,1]\to\Omega$ between any two points $\gamma(0)=x\in S$ and $\gamma(1)=y\in\Omega\setminus S$. Now take $r_0 := \inf_{0\leq t\leq 1} \dist(\gamma(t),\de\Omega)$. Note that since $\gamma([0,1])$ is closed, we have $r_0>0$. Consider $T = \sup\{0\leq t\leq 1:\, \gamma(t) \in S\}$ and note that $T<1$ since $\gamma(1)\not\in S$. But since we know that $B_{\delta_0 r_0}(\gamma(T)) \subset B_{\delta_0\dist(\gamma(T),\de\Omega)}(\gamma(T)) \subset S$ and by continuity  we can verify that for small enough $\epsilon$ we have $\gamma(T+\epsilon) \in B_{\delta_0r_0}(\gamma(T))$, so that $f\equiv Q\a{0}$ in a neighborhood of $\gamma(T+\epsilon)$ and therefore $\gamma(T+\epsilon)\in S$. This is in contradiction with the definition of $T$ and our claim follows. Hence $S=\Omega$ and $f=Q\llbracket0\rrbracket$ in $\Omega$.

        Notice that \cref{thm:strong} follows since the assumption therein guarantees that $f=Q\a{0}$ and moreover that $\kappa_{0,f}=\infty$.
    \end{proof}
\end{lemma}

Having ruled out the case $\kappa_x=\infty$, we can show that for a nontrivial Dir-stationary multivalued function $\kappa_x$ is well defined.

\begin{lemma}\label{vanihsing-order-limit-lemma}
    Let $f\in W^{1,2}(\Omega, \Iq(\R^m))$ be a nontrivial, Dir-stationary multivalued function. Then for every $x\in \Omega$, the vanishing order $\kappa_x$ is well defined as 
    \[
    \kappa_x=\lim_{r\to 0} \frac{\log\left(\fint_{B_{2r}\setminus B_r(x)} |f|^2 \right)}{2\log(r)}
    \]
    
    \begin{proof}
        For any point $x$ take a sequence of radii $\{r_j\}_{j=1}^\infty$ realizing the lim-sup $\kappa_x<\infty$. Taking $j$ large enough so that
        \begin{align*}
            \frac{\log\left(\fint_{B_{2r_j}\setminus B_{r_j}(x)} |f|^2 \right)}{2\log(r_j)} \geq \kappa_x - \frac{\epsilon}{2} 
        \end{align*}
        we see, since $\log(r_j)\leq 0$, that
        \[
        \fint_{B_{2r_j}\setminus B_{r_j}(x)} |f|^2 \leq r_j^{2\kappa_x - \epsilon } \,.
        \]
        Let $r>0$ be such that $r_{j-1} \geq r \geq r_j$. Then the three sphere inequality for $\tau = \kappa_x + \frac{n}2 - \epsilon$ states
        \begin{align*}
            \frac{\fint_{B_{2r}\setminus B_{r}(x)} |f|^2}{(\log(r)^2+1)r^{2\kappa_x-2\epsilon}} \leq C(\kappa_x)\left( \frac{\fint_{B_{2r_j}\setminus B_{r_j}(x)} |f|^2}{r_j^{2\kappa_x-2\epsilon}} + \frac{\fint_{B_{2r_{j-1}}\setminus B_{r_{j-1}}(x)} |f|^2}{r_{j-1}^{2\kappa_x-2\epsilon}}\right)\,.
        \end{align*}
        Note that since $\kappa_x$ is bounded, then one can fix $\tau_0$, to make the constant in the three sphere inequality uniformly bounded throughout our argument here. In the rest of the proof, we drop the dependence of $C(\kappa_x)$. Hence we see that for all $j > j_0(\epsilon)$ large enough and $r$ as above:
        \begin{align*}
            \frac{\fint_{B_{2r}\setminus B_{r}(x)} |f|^2}{(\log(r)^2+1)r^{2\kappa_x-2\epsilon}} \leq C \,r_{j-1}^{\epsilon}\leq C\, r^{\epsilon}\,.
        \end{align*}
        Taking logarithm, for $r$ sufficiently small we have
        \begin{align*}
            \frac{\log\left(\fint_{B_{2r}\setminus B_r(x)}|f|^2\right)}{2\log(r)} \geq \kappa_x - \epsilon - C\frac{\log(\log(r))}{\log(r)} \geq \kappa_x - 2\epsilon\,.
        \end{align*}
        This means that 
        \begin{align*}
            \kappa_x = \lim_{r\to0}\frac{\log\left(\fint_{B_{2r}\setminus B_r(x)}|f|^2\right)}{2\log(r)} \,.
        \end{align*}
        and we are done.
    \end{proof}
\end{lemma}

\begin{remark}[Frequency and vanishing order coincide]
    As a consequence of \cite[EQ (3.42)]{DS} and \cref{vanihsing-order-limit-lemma}, for Dir-stationary multi-valued functions, the vanishing order and the frequency (defined via a linear cut-off as in \cite[Definition 3.1]{DLS_Blowup}) coincide:
    \begin{align}
        \kappa_{x,f} = I_{x,f}\,.
    \end{align}
    Note that the same conclusion of \cite[EQ (3.42)]{DS} applies, hence integrating from $s$ to $t$ yields:
    \begin{align*}
        \log\left(\frac{H(r)}{r^{n-1}}\right) - \log\left(\frac{H(s)}{s^{n-1}}\right) = \int_{s}^{r} \frac{2I(\tau)}{\tau} d\tau 
    \end{align*}
    For $s\leq r\leq \delta$ small enough, we know that the frequency is saturated $|I(r)-I_{x,f}|\leq\epsilon$, hence
    \begin{align*}
         \log\left(\frac{H(r)}{r^{n-1}}\right) - \log\left(\frac{H(s)}{s^{n-1}}\right) = (2I_{x,f} + O(\epsilon)) \log\left(\frac{r}{s}\right)\,.
    \end{align*}
    By the definition of vanishing order, for small enough $\delta>0$:
    \begin{align*}
         \log\left(\frac{H(r)}{r^{n-1}}\right) - \log\left(\frac{H(s)}{s^{n-1}}\right) =(2\kappa_{x,f}+O(\epsilon))\log\left(\frac{r}{s}\right)\,.
    \end{align*}
    This shows that:
    \begin{align*}
        |I_{x,f} - \kappa_{x,f}| = O(\epsilon)\,,
    \end{align*}
    for all $\epsilon>0$ and we conclude.
\end{remark}

\subsection{Properties of the vanishing order}

In this section we will show upper-semi continuity of $\kappa_x$, and indeed that $\kappa_x$ is uniformly bounded on $B_1$, and also study the homogeneity of subsequential blow-up limits.

\begin{lemma}\label{lemma:upper-semi-cont}
    Let $f$ be a Dir-stationary multi-valued functions with zero average, then vanishing degree $x \mapsto \kappa_x$ is upper-semi continuous.
\end{lemma}
    
    \begin{proof}
        We can compare close-by points by using the three-sphere inequality. Without loss of generality via a translation we only need to prove upper-semi continuity for the origin, that is we aim to show that for every $\epsilon>0$ there is $\delta>0$ such that $\kappa_0 \geq \kappa_x - \epsilon$ for any $x\in B_\delta$. This will follow if we can show that for any $x\in B_\delta$ there exists $r_0(x,\epsilon)>0$ such that $\forall r\leq r_0(x,\epsilon)$ we have
        \begin{align*}
            \int_{B_{2r}\setminus B_r (x)} |f|^2 \geq r^{2\kappa_0 + n + 2\epsilon}\,.
        \end{align*}
        
        To achieve this, note that for any three radii $r \ll K_1\delta \leq K_2\delta$ with large $2< K_1 < K_2/2$ (to be chosen later), \cref{three-sphere-multivalued-corollary} (noting uniform bounds on $\tau$ depending on $\kappa_0$ throughout the argument) with $\tau = \kappa_0 + \frac{n}{2} + \epsilon$ says:
        \begin{align}\label{upper-semi-three-sphere-temp}
            \frac{\fint_{B_{K_1\delta}\setminus B_{\frac{K_1\delta}{2}}(x)} |f|^2 }{(\log(K_1/K_2)^2 + 1)(K_1\delta)^{2\kappa_0+2\epsilon}}  - C \frac{\fint_{B_{K_2\delta}\setminus B_{\frac{K_2\delta}{2}}(x)} |f|^2 }{(K_2\delta)^{2\kappa_0+2\epsilon}}\leq C \frac{\fint_{B_{2r}\setminus B_{r}(x)} |f|^2 }{r^{2\kappa_0+2\epsilon}}\,.
        \end{align}
        Note that here we take $r$ so small so that $\log(K_2/K_1) \leq \log(K_1\delta/r)$. Now for $x \in B_\delta$ and large enough $K_1$ we have that $B_{\frac45 K_1\delta}\setminus B_{\frac35K_1\delta}(0) \subset B_{K_1\delta}\setminus B_{\frac{K_1\delta}{2}}(x)$, hence:
        \begin{align}\label{upper-semi-temp-1}
            \fint_{B_{K_1\delta}\setminus B_{\frac{K_1\delta}{2}}(x)} |f|^2 \geq \fint_{B_{\frac45 K_1\delta}\setminus B_{\frac35K_1\delta}(0)} |f|^2 \geq C(K_1\delta)^{2\kappa_0+\epsilon}\,,
        \end{align}
        for small enough $\delta>0$. This follows from the vanishing order definition. Similarly, for large enough $K_2$ and small enough $\delta$ we also have that $B_{\frac45K_2\delta}\setminus B_{\frac35K_2\delta} (x) \subset B_{K_2\delta}\setminus B_{\frac{K_2\delta}{2}} (0)$, hence we see that:
        \begin{align}\label{upper-semi-temp-2}
            \fint_{B_{\frac45K_2\delta}\setminus B_{\frac35K_2\delta} (x)} |f|^2 \leq \fint_{B_{K_2\delta}\setminus B_{\frac{K_2\delta}{2}} (0)} |f|^2 \leq C(K_2\delta)^{2\kappa_0-\epsilon}\,.
        \end{align}
        Putting together \cref{upper-semi-temp-1,upper-semi-temp-2} with the three sphere inequality \cref{upper-semi-three-sphere-temp}, we conclude that:
        \begin{align*}
            C\frac{\fint_{B_{2r}\setminus B_{r}(x)} |f|^2 }{r^{2\kappa_0+2\epsilon}} \geq C\left(\frac{(K_1\delta)^{-\epsilon}}{\log(K_2/K_1)^2}-(K_2\delta)^{-3\epsilon}\right)\,.
        \end{align*}
        We can take $K_2$ large enough such that $\frac{(K_1\delta)^{-\epsilon}}{\log(K_2/K_1)^2}>(K_2\delta)^{-3\epsilon}$ and we conclude that:
        \begin{align*}
            C \frac{\fint_{B_{2r}\setminus B_{r}(x)} |f|^2 }{r^{2\kappa_0+2\epsilon}} \geq C(K_1,K_2,\delta,\epsilon)\,,
        \end{align*}
        for all $r>0$ sufficiently small. Since the right hand side is independent of $r$, 
        we take a logarithm, and conclude that:
        \begin{align*}
            \kappa_x \leq \kappa_0 + \epsilon\,.
        \end{align*}
        This is indeed the desired conclusion. We see in fact that for any $x\in B_1$:
        \begin{align*}
            \kappa_x \geq \limsup_{y\to x} \kappa_y\,.
        \end{align*}

    \end{proof}

The upper semi-continuity of \cref{lemma:upper-semi-cont} together with strong unique continuation in \cref{lemma:strong-unique-continuation} imply that the vanishing order $\kappa_x$ is uniformly bounded in $x$.

In the next proposition, we prove that around any point $x\in B_1$ there exists a sequence of radii such that the solution becomes $\kappa_x$ homogeneous along the sequence. 

\begin{theorem}\label{almost-homogen-thm}
    Let $f\in W^{1,2}(B_1;\A_Q)$ be a Dir-stationary multi-valued function. Then around any point $x$ there exists a sequence of radii $r_j\to 0$ and vanishing constants $\epsilon_j\to0$ such that:
    \begin{align}\label{eq:homoimpr}
        \int_{B_{2r_j}\setminus B_{r_j}(x)} \sum_{i=1}^Q |Df_i.x - \kappa_x f_i|^2 \leq \epsilon_j \int_{B_{2r_j}\setminus B_{r_j}(x)} \sum_{i=1}^Q |f_i|^2\,.
    \end{align}
    \end{theorem}
    
    \begin{proof}
        It is enough to prove the result for the origin. The proof is by contradiction. Indeed assume that there exists $\epsilon_0>0$ such that for all small radii $r \leq r_0$ we have:
        \begin{align}\label{contradict-1}
            \int_{B_{2r}\setminus B_r(x)} \sum_{i=1}^Q |Df_i.x - \kappa_x f_i|^2\geq \epsilon_0 \int_{B_{2r}\setminus B_r(x)} \sum_{i=1}^{Q}|f_i|^2\,.
        \end{align}
        The idea is from \cite{KRS-carleman} and is as follows.  Instead of the precise weights used in \cref{thm:Carleman-proposition}, we use the following
        \begin{align*}
            \text{inner variations: }xe^{-2\tau \phi(\log(|x|))} \text{ \& outer variation: }ue^{-2\tau \phi(\log(|x|))}\,.
        \end{align*}
        Take a compactly supported function $\chi \in C^\infty_c(B_1\setminus\{0\})$; The inner variation \cref{eq:inner} implies:
        \begin{align}
        \begin{aligned}\label{homogen-1}
            0 &= \mathcal{I}\left(f,\chi xe^{-2\tau\phi(\log(|x|))}\right) \\ &= \int_{B_1} \chi\frac{n-2\tau\phi'-2}{2}\sum_{i=1}^Q\left(|Df_i|^2 + 2\tau|\de_r f_i|^2\right) e^{-2\tau\phi(\log(|x|))}\\
            & + \int_{B_1} \sum_{i=1}^Q\left((D\chi.x) \frac{|Df_i|^2}{2} - (D\chi.Df_i)(Df_i.x)\right)e^{-2\tau\phi(\log(|x|))}\,.
            \end{aligned}
        \end{align}
        For the outer variation \cref{eq:outer}, with the same smooth cut-off $\chi$ we see that:
        \begin{align}
        \begin{aligned}\label{homogen-2}
            0&=\mathcal{O}(f,\chi u e^{-2\tau\phi(\log(|x|))}) \\&= \int_{B_1} \chi \sum_{i=1}^Q \left(|Df_i|^2 - 2\tau\phi'\frac{\de_rf_i f_i}{|x|}\right) e^{-2\tau\phi(\log(|x|))}\\
            &+\int_{B_1} \sum_{i=1}^Q (D\chi.Df_i)f_i e^{-2\tau\phi(\log(|x|))}\,.
            \end{aligned}
        \end{align}
        Then we multiply \cref{homogen-2} by $\eta=\frac{2\tau-n+2}{2}$ and add to \cref{homogen-2} to get the following inequality:
        \begin{align*}
            &\int_{B_1} \chi\sum_{i=1}^Q \left(|\de_rf_i|^2 - \eta^2 \frac{|f_i|^2}{|x|^2}\right) e^{-2\tau\phi(\log(|x|))} \\ \leq C &\int_{B_1}|D\chi|\sum_{i=1}^Q\left(|x||Df_i|^2 + \frac{|f_i|^2}{|x|}\right)e^{-2\tau\phi(\log(|x|))} \\ &+ C \|1-\phi'\|_\infty\int_{B_1}\chi\sum_{i=1}^Q\left(|Df_i|^2 + \frac{|f_i|^2}{|x|^2}\right)e^{-2\tau\phi(\log(|x|))}\,.
        \end{align*}
        Performing the same integration by parts in the proof of \cref{thm:Carleman-proposition} in \cref{temp-1} we see that:
        \begin{align}\label{modified-carleman}
        \begin{aligned}
            &\int_{B_1} \chi\sum_{i=1}^Q \left|\de_rf_i-\frac{f_i}{|x|}\right|^2e^{-2\tau\phi(\log(|x|))} \\ \leq C &\int_{B_1}|D\chi|\sum_{i=1}^Q\left(|x||Df_i|^2 + \frac{|f_i|^2}{|x|}\right)e^{-2\tau\phi(\log(|x|))} \\+ C&\left( \|1-\phi'\|_\infty + \|\phi''\|_\infty\right)\int_{B_1}\chi\sum_{i=1}^Q\left(|Df_i|^2 + \frac{|f_i|^2}{|x|^2}\right)e^{-2\tau\phi(\log(|x|))}\,.
            \end{aligned}
        \end{align}
        Now we use the contradiction assumption (\cref{contradict-1}). The idea is to \textit{bend} slightly the graph of $\phi(t)$ with a convexity of size $\epsilon_0>0$. With that we gain a three-sphere inequality with $\epsilon_0>0$ more weight for the left hand side. Using this we gain a contradiction with the fact that the vanishing order is a limit in \cref{vanihsing-order-limit-lemma}. Now for any two radii $r_1 \ll r_2 \leq r_0(x)$ we introduce $\phi_\delta$ as follows:
        \begin{align*}
                \begin{cases}
                    \phi_\delta(t) \geq (1-\delta)t &\text{for } \log(r_1)\leq t\leq \log(2r_1) \text{ or } \log(r_2)\leq t\leq \log(2r_2)\,,\\
                    \phi_\delta(t) \leq (1+2\delta)t &\text{for } \log(\sqrt{r_1r_2}) \leq t \leq \log(2\sqrt{r_1r_2})\,,\\
                    |\phi'| + |\phi''| \leq C\delta\,.
                \end{cases}
        \end{align*}
        Moreover we put:
        \begin{align*}
            \tau= \frac{2\kappa_x + n - 2 }{2} \Rightarrow \eta=\frac{2\tau - n + 2}{2} = \kappa_x
        \end{align*}
        Now using \cref{modified-carleman} and \cref{contradict-1} and the Cacciopoli inequality on dyadic annuli between $r_1$ and $r_2$ we see that:
        \begin{align*}
            \int_{B_1} \chi\frac{|f|^2}{|x|^2}e^{-2\tau\phi_\delta(\log(|x|))}  \leq \frac{C}{\epsilon_0-C\delta} \int_{B_1}(|D\chi| + |D^2\chi|)\frac{|f|^2}{|x|^2}e^{-2\tau\phi_\delta(\log(|x|))}
        \end{align*}
        Now take $\chi$ to be the smooth cut-off such that $\chi=1$ on $B_{r_2}\setminus B_{2r_1}(x)$ and it linearly decreases to $0$ on $B_{r_1}(x)$ and $\left(B_{2r_2}(x)\right)^c$. Then we can see that for $\delta = c\epsilon_0$ for small enough $c>0$, there exists a constant $C(\epsilon_0,x)$ such that:
        \begin{align*}
            \frac{\int_{B_{2\sqrt{r_1r_2}}\setminus B_{\sqrt{r_1r_2}}(x)} |f|^2 }{\sqrt{r_1r_2}^{2\kappa_x + 2c\epsilon_0+ n}} \leq C(\epsilon_0,x) \left[\frac{\int_{B_{2r_1}\setminus B_{r_1}(x)}|f|^2}{r_1^{2\kappa_x-c\epsilon_0+n}}+\frac{\int_{B_{2r_2}\setminus B_{r_2}(x)}|f|^2}{r_2^{2\kappa_x-c\epsilon_0+n}}\right]\,.
        \end{align*}
        By the definition of vanishing order $\kappa_x$, for any $\eta>0$ there exists $r_0(\eta,x)$ such that for all radii $r\leq r_0(\eta,x)$ we have:
        \begin{align*}
            r^{\kappa_x + \eta} \leq \left(\fint_{B_{2r}\setminus B_r(x)}|f|^2\right)^{\frac12} \leq r^{\kappa_x - \eta}\,.
        \end{align*}
        Combining the last two displays, we arrive at:
        \begin{align*}
            \sqrt{r_1r_2}^{2\eta-2c\epsilon_0} \leq C(\epsilon_0)\left(r_1^{c\epsilon_0-2\eta} + r_2^{c\epsilon_0-2\eta}\right)\,.
        \end{align*}
        We can take $\eta$ small enough so that $2\eta - 2c\epsilon_0 < 0$ and $c\epsilon_0-2\eta > 0$, hence the right hand side becomes bounded by $C(\epsilon_0)$. Then taking $r_1\ll r_2$ small enough we reach a contradiction and conclude.
    \end{proof}

\begin{remark}
    The previous proof actually implies that for every sequence of radii $r_j\to 0$ there exists a subsequence $r_{j_k}\to 0$ and a sequence of constants $\eps_{j_k}\to 0$ for which \eqref{eq:homoimpr} holds. 
\end{remark}

\section{Proof of the dimensional bound: \cref{thm:dimension}}

The dimensional bound in \cref{thm:dimension} follows as in \cite[Proof of Theorem 0.11]{DS} from the following Lemma.

\begin{lemma}\label{lem:dimred} Let $\Omega$ be connected and $f \in W^{1,2}(\Omega, \Iq(\R^m))$ be Dir-minimizing. Then, either $f= Q \a{\xi}$ with $\xi \colon \Omega \to \R^n$ harmonic in $\Omega$, or the set 
\[
\Sigma_{Q,f}:=\{x\in \Omega\,:\,f(x)=Q\a{y},\,y\in \R^n\}
\](which is relatively closed in $\Omega$ since $f$ is continuous) has Hausdorﬀ dimension at most $n-2$ and it is locally
finite for $n=2$.
\end{lemma}

First we notice that we can assume $\xi, y\equiv 0$ in \cref{lem:dimred}, since we can subtract the average of $f$ by \cite[Lemma 3.23]{DS}. Next we define the blow-up sequence
\[
f_{y,\rho}(x):=\frac{\rho^{\frac{n}2}f(y+\rho x)}{\sqrt{\int_{B_\rho(y)}|f|^2}}
\] 
The proof of \cref{lem:dimred} will then follow as in the proof of \cite[Proposition 3.22]{DS} replacing Theorem 3.19 therein with the following

\begin{lemma}\label{lem:homobu}
Let $f\in W^{1,2}(B_1, \Iq(\R^m))$ be Dir-minimizing. Assume  $f(0)=Q\a{0}$ and $\|f\|_{L^2(B_\rho))}>0$ for every $\rho\leq 1$.  Then, for any sequence $(f_{\rho_k})_k$, with $\rho_k \downarrow 0$, a subsequence, not relabeled, converges locally uniformly to a function $g\colon \R^n \to \Iq(\R^m)$ with the following
properties:
\begin{enumerate}
\item $\|g\|_{L^2(B_1)}=1$ and $g|_\Omega$ is Dir-minimizing for any bounded $\Omega$;
\item  $g(x) = |x|^\alpha g\left(\frac{x}{|x|}\right)$, where $\alpha= k_{0,f} >0$ is the vanishing order of $f$ at $0$.
\end{enumerate}
\end{lemma}

Before we prove \cref{lem:homobu}, we need the following doubling estimate to guarantee that the sequence $f_{\rho}$ is uniformly bounded.

\begin{proposition}[Doubling estimate]\label{Doubling-prop}
    Let $f\in W^{1,2}(B_1^n,\A_Q)$ be a Dir-stationary multifunction. Then for any point $x\in B_1$ there exists a radius $r_x>0$ and a constant $C_x$ such that for all radii $r \leq r_x$ we have:
    \begin{align*}
        \int_{B_{2r}(x)} |f|^2 \leq C_x \int_{B_r(x)}|f|^2\,.
    \end{align*}
  \end{proposition}

    \begin{proof}
        It is enough to prove that:
        \begin{align*}
            \int_{B_{2r}\setminus B_r(x)} |f|^2 \leq C_x \int_{B_r(x)} |f|^2\,.
        \end{align*}
        The strategy is to use the three sphere inequality for the three radii $\epsilon \leq 2\epsilon \ll r_3$ for $r_3>0$ to be chosen later:
        \begin{align*}
            C\frac{\int_{B_{2\epsilon}\setminus B_\epsilon(x)} |f|^2}{(2\epsilon)^{2\tau}} \leq \frac{\int_{B_{\epsilon}\setminus B_{\epsilon/2}(x)}|f|^2}{\epsilon^{2\tau}} + \frac{\int_{B_{2r_3}\setminus B_{r_3}(x)}|f|^2}{r_3^{2\tau}} 
        \end{align*}
        Multiplying and rearranging we see that:
        \begin{align}\label{pre-doubling}
            \int_{B_{2\epsilon}\setminus B_\epsilon(x)} |f|^2 \leq C2^{2\tau}\int_{B_{\epsilon}\setminus B_{\epsilon/2}(x)}|f|^2 + C \left(\frac{\epsilon}{r_3}\right)^{2\tau} \int_{B_{2r_3}\setminus B_{r_3}(x)} |f|^2\,.
        \end{align}
        Now by \cref{vanihsing-order-limit-lemma}, we know that for any $\eta>0$ there exists a radius $r_0(\eta,x)$ small enough such that for all radii $r\leq r_0(\eta,x)$ we have that:
        \begin{align*}
            r^{2\kappa_x + n + \eta} \leq \int_{B_{2r}\setminus B_r(x)} |f|^2 \leq r^{2\kappa_x + n- \eta}\,.
        \end{align*}
        We take $r_3 \leq r_0(\eta,x)$, and we see that:
        \begin{align*}
            C \left(\frac{\epsilon}{r_3}\right)^{2\tau} \int_{B_{2r_3}\setminus B_{r_3}(x)} |f|^2 \leq \left(\frac{\epsilon}{r_3}\right)^{2\tau} r_3^{2\kappa_x+n-\eta} 
        \end{align*}
        Take $\tau=\frac{2\kappa_x+n+ 4\eta}{2}$ and we see that the above display reads:
        \begin{align*}
            C \left(\frac{\epsilon}{r_3}\right)^{2\tau} \int_{B_{2r_3}\setminus B_{r_3}(x)} |f|^2 \leq \epsilon^{2\kappa_x+n+3\eta} \leq \epsilon^{2\eta}\int_{B_{2\epsilon}\setminus B_\epsilon(x)}|f|^2\,.
        \end{align*}
        Taking $\epsilon>0$ small enough we see that we can we can reabsorb the second term on the right hand side of \cref{pre-doubling} and end up with the desired estimate:
        \begin{align*}
            \int_{B_{2\epsilon}\setminus B_\epsilon(x)} |f|^2 \leq C(x,n)\int_{B_{\epsilon}(x)}|f|^2\,.
        \end{align*}
    \end{proof}

\begin{proof}[Proof of \cref{lem:homobu}] We consider any ball $B_N$ of radius $N$ centered at $0$. It follows from \cref{Doubling-prop} and the Caccioppoli inequality \cref{eq:caccio} that $\|f_\rho\|_{W^{1,2}(B_N)}$ is uniformly bounded in $\rho$. Hence, the functions $f_\rho$ are all Dir-minimizing and \cite[Theorem 3.9]{DS} implies that they are locally equi-H\"older continuous. Since $f(0) = Q \a{0}$, the $f_\rho$’s are also locally uniformly bounded and the
Ascoli–Arzel\'a theorem yields a subsequence (not relabeled) converging uniformly on compact subsets of $\R^n$ to a continuous $Q$-valued function $g$. This implies easily the weak convergence (see \cite[Definition 2.9]{DS}), so we can apply \cite[Proposition 3.20]{DS} and conclude
(1) (note that $\|f_\rho\|_{L^2(B_1)} = 1$ for every $\rho$). 

Next notice that up to choosing a further subsequence $g$ is $k_{0,f} (0)$-homogeneous by \cref{lem:homobu}. Finally, assume by contradiction that $k_{0,f} (0) = 0$. Then, by what shown so far, the blowups converge to a continuous $0$-homogeneous function $g$, with $g(0) = Q\a{0}$. This implies that $g \equiv Q \a{0}$, a contradiction to $\|g\|_{L^2(B_1)}=1$.
\end{proof}

\section{The 2-dimensional case}
For the reader's convenience we give a proof of the $2$-dimensional case combining the monotonicity of the Weiss energy with the reasoning above on the vanishing order of a Dir-minimizing multivalued map. We start with the following:

\begin{proposition}[Weiss' monotonicity formula]\label{p:Weiss}
If $f\in W^{1,2}(\Omega,\Iq(\R^m))$ and $B_r (x)\subset \Omega$, then we define the \emph{Weiss Boundary adjusted energy} by
\[
W_{x,f}(r)=W(f_{x,r}):=\frac{1}{r^{m+2\kappa_{x,f}-2}}\int_{B_r(x)}|Df|^2-\frac{\kappa_{x,f}}{r^{m+2\kappa_{x,f}-1}}\int_{\partial B_r(x)}|f|^2\,.
\]
Then the map $r\mapsto W_{x,f} (r)$ is absolutely continuous and
\begin{equation}\label{e:Weiss_monot}
\frac{d}{dr} W_{x,f} (r) =\frac{m+2\kappa_{x,f}-2}{r}(W(f^{\kappa_{x,f}}_{x,r})-W(f_{x,r}))+\frac1r\int_{\partial B_1}\sum_{i=1}^Q\left|(Df_{x,r})_i\cdot x-\kappa_{x,f}(f_{x,r})_i\right|^2\, ,
\end{equation}
where $f^\kappa_{x,r}(y):=|y|^\kappa\,f_{x,r}(x/|x|)$ is the $\kappa$-homogeneous extension of the trace of $f_{x,r}$ in $B_1$. 
In particular if $f$ is Dir-minimizing in $\Omega$ then 
\begin{enumerate}
    \item $\frac{d}{dr} W_{x,f} (r)\geq 0$ and so there exists $W_{x,f}(r)\geq W_{x,f}(0)=\lim_{r\downarrow 0}W_{x,f}(r)= 0$;
    \item $\frac{d}{dr} W_{x,f} (r)\equiv 0$ if and only if $f_{x,r}$ is $\kappa_{x,f}$-homogeneous.
\end{enumerate}
\end{proposition}

\begin{proof} The proof of the above is standard and can be found for instance in \cite[Section 9]{Ve} (adjusting the constants therein to account for $\kappa$). In particular notice $W_{x,f}(0)=0$ follows from \cref{lem:homobu} since $W_{x,f}(0)=W_{0,g}(1)$, where $g$ is a $\kappa_{x,f}$ homogeneous function and so $W_{0,g}(1)=0$.
\end{proof}

Next we have the following standard epiperimetric inequality.

\begin{lemma}[Epiperimetric Inequality]\label{l:epiperimetric} There is a positive constant $\delta$, depending only on $Q$, with the following property. Assume $f\in W^{1,2} (B_1, \Iq(\R^m))$, $B_1\subset \R^2$, is $\Dir$-minimizing. Then, for $f^\kappa$, the $\kappa$ homogeneous extension of $f$, we have
\begin{equation}\label{e:epiperimetric}
\int_{B_1} \left(|Df^\kappa|^2 - |Df|^2\right) \geq \delta\, W_{0,f} (1)\, . 
\end{equation} 
\end{lemma}

\begin{proof} 
    In short, we construct competitors by unwinding the boundary conditions, extending harmonically inside and rewinding again. We use the same notation as in \cite[Proposition 5.2]{DS}; precisely, fix a radius $r$ and let $f(re^{i\theta}) = g(\theta) = \sum_{j=1}^{J}\llbracket g_j(\theta)\rrbracket$ be an irreducible decomposition as in \cite[Proposition 1.5]{DS}. Then for each $g_j$ we can find $\gamma_j:S^1\to\R^n$ such that:
    \begin{align*}
        g_j(\theta) = \sum_{i=1}^{Q_j} \left\llbracket \gamma_i\left(\frac{\theta + 2\pi i}{Q_j}\right)\right\rrbracket\,.
    \end{align*}
    Now take the Fourier decomposition of $\gamma_j$:
    \begin{align*}
        \gamma_j(\theta) = \frac{a_{j,0}}{2} + \sum_{\ell = 1}^{\infty} \left[ a_{j,\ell} \sin(\ell\theta) + b_{j,\ell} \cos(\ell\theta) \right]\,,
    \end{align*}
    and its harmonic extension:
    \begin{align*}
        \zeta_j(\rho,\theta) = \frac{a_{j,0}}{2} + \sum_{\ell=1}^{\infty} \rho^\ell\left[ a_{j,\ell} \sin(\ell\theta) + b_{j,\ell} \cos(\ell\theta) \right]\,.
    \end{align*}
    Calculating the Dirichlet energy as in \cite[EQ (5.18)]{DS}, we obtain:
    \begin{align*}
        \int_{B_r}|Df|^2 \leq \sum_j \textrm{Dir}(\zeta_j,B_r) = \pi \sum_{j}\sum_\ell \ell r^{2\ell}\left(|a_{j,\ell}|^2 + |b_{j,\ell}|^2\right)\,. 
    \end{align*}
    We also calculate the Dirichlet energy for the $\kappa$-homogeneous extension:
    \begin{align*}
        \int_{B_r}|Df^\kappa|^2 = \sum_j \textrm{Dir}(\gamma^{\kappa}_j,B_r) = \pi r^{2\kappa}\sum_{j}\sum_\ell \left[\frac{\kappa}{2} + \frac{\ell^2}{2\kappa}\right]\left(|a_{j,\ell}|^2 + |b_{j,\ell}|^2\right)\,.
    \end{align*}
    Hence we can estimate for $r=1$:
    \begin{align}\label{Weiss-lower-bound}
        \int_{B_1} |Df^\kappa|^2 - |Df|^2 \geq
        \pi\sum_j\sum_\ell \left[ \frac{\kappa}{2} + \frac{\ell^2}{2\kappa} - \ell \right]\left(|a_{j,\ell}|^2 + |b_{j,\ell}|^2\right)\,.
    \end{align}
    Now we calculate the Weiss energy using competitors:
    \begin{align}\label{Weiss-upper-bound}
        W(f_{x,1})  \leq \sum_j \textrm{Dir}(\zeta_j,B_1) - \kappa \int_{\de B_1} |f|^2 \leq \pi \sum_j \sum_{\ell} \left[\ell - \kappa Q_j\right]\left(|a_{j,\ell}|^2 + |b_{j,\ell}|^2\right)\,.
    \end{align}
    Since $Q_j\geq 1$, it is enough to show that there exists some $\delta$ such that:
    \begin{align}\label{epi-temp-1}
        \frac{\ell^2}{2\kappa} + \frac\kappa2 - \ell \geq \delta(\ell - \kappa )\,.
    \end{align}
    It is easy to verify that the following choice for $\delta$ satisfies \cref{epi-temp-1}
    \begin{align*}
        \delta = \frac{\lfloor\kappa\rfloor + 1 -\kappa}{2\kappa}
    \end{align*}
    Putting together \cref{Weiss-lower-bound,Weiss-upper-bound}, we conclude that:
    \begin{align*}
        \int_{B_1} |Df^\kappa|^2 - |Df|^2 \geq \delta W_{0,f}(1)\,.
    \end{align*}
\end{proof}

Combining the two propositions above we have the following

\begin{lemma}[Uniqueness of tangent map]\label{lem:unique} Let $f \in W^{1,2} (B_1, \Iq(\R^m))$ be a Dir-minimizing Q-valued functions, with $\Dir(f,B_1) > 0$ and $f(0) = Q\a{0}$. Then, the maps $f_{x,r}$ converge locally uniformly to a unique tangent mat $g$. 
\end{lemma}

\begin{proof} The proof follows from a standard reasoning, see for instance \cite[Lemma 12.14 \& Proposition 2.14]{Ve}.
\end{proof}

The proof of the $2$-dimensional case of \cref{thm:dimension} then follows in the same way as in \cite[Subsection 5.3]{DS}, replacing Theorem 5.3 therein with \cref{lem:unique}.

\bibliographystyle{acm}
\bibliography{references.bib}

\end{document}